\documentclass[10pt]{article} 
\usepackage[spanish]{babel}
\usepackage{amsmath,amssymb,amsfonts}
\usepackage{hyperref} 
\hypersetup{hidelinks,colorlinks,breaklinks=true,urlcolor=blue,citecolor=green,linkcolor=green,bookmarksopen=false}

\title{Sobre \emph{Introducci\'on al An\'alisis Matem\'atico} de Mario O. Gonz\'alez}
\author{Luis Giraldo Gonz\'alez Ricardo} 

\begin{document}

\maketitle 
\begin{abstract}
	En el presente art\'iculo, tomamos como pretexto los 75 a\~nos de la publicaci\'on de \emph{Introducci\'on al An\'alisis Matem\'atico} de Mario O. Gonz\'alez para estudiar su trascendencia y elementos notables. Adem\'as, se comentan algunas cuestiones novedosas de la biograf\'ia del autor antes de la publicaci\'on de ese libro.
	
	In this article, we use the 75th anniversary of the publication of \emph{Introducci\'on al An\'alisis Matem\'atico} by Mario O. Gonz\'alez to study its importance. Besides, we give new biographical data from his author in the previous years to the publication of the book.	
\end{abstract}

\thispagestyle{empty} 


\selectlanguage{spanish}
\section*{Introducci\'on} 
La etapa colonial cubana fue adversa para el desarrollo de la cultura matem\'atica. Por una parte, la mayor\'ia de la poblaci\'on cubana era analfabeta o pose\'ia escasa instrucci\'on; y por otra, los pocos que ten\'ian acceso a la Universidad de La Habana –\'unica de su tipo en el pa\'is–, encontraban que el nivel de las asignaturas relacionadas con la matem\'atica era muy bajo. En este sentido puede consultarse \cite{CarlosConcha}. Con el nacimiento de la Rep\'ublica, gracias al Plan Varona, se movieron nuevos aires en la educaci\'on superior cubana. A pesar de ello, la situaci\'on de la «reina de las ciencias» no cambi\'o en demas\'ia. Baste se\~nalar que el catedr\'atico de An\'alisis de la Facultad de Letras y Ciencias de la Universidad de La Habana, el Ing. Jos\'e A. Villal\'on (1864-193?), ten\'ia numerosas deficiencias en lo que a sus conocimientos de an\'alisis matem\'atico respecta. Fue el profesor de origen catal\'an Claudio Mim\'o y Caba (1844-1929), quien aglutin\'o a su alrededor a los primeros doctores en Ciencias F\'isico-Matem\'aticas de la naci\'on \cite{CarlosConcha}. 

Fue el inicio de la actividad docente de Pablo Miquel y Merino (1877-1944), en el curso 1912-13, lo que marc\'o un cambio cualitativo en la instrucci\'on del an\'alisis en la Universidad de La Habana, lo que es visible en su libro \emph{Elementos de \'algebra superior} (1914). All\'i, se tratan los conceptos b\'asicos del an\'alisis de manera did\'actica, y con el nivel de rigor t\'ipico de la segunda d\'ecada del siglo \textsc{xx}, superior en muchos sentidos a los otros textos empleados en la ense\~nanza del c\'alculo \cite{ConchaCarlos}. Gracias al arduo trabajo del profesor Miquel, la cultura matem\'atica cubana fue en ascenso con el transcurso del siglo. Bajo su auspicio, se formaron muchos de los pedagogos que dieron prestigio los cursos de matem\'atica en los Institutos de Segunda Ense\~nanza. Sin duda, el m\'as reconocido de sus pupilos fue Mario O. Gonz\'alez (1913-1999), quien se convirti\'o en el matem\'atico cubano m\'as destacado a nivel internacional antes del triunfo revolucionario de 1959. Por mucho tiempo los textos de Mario O. Gonz\'alez fueron utilizados en la ense\~nanza secundaria y preuniversitaria de Cuba –por lo que su influencia en la formaci\'on matem\'atica de los cubanos trascendi\'o su presencia en el pa\'is– y entre ellos, los m\'as conocidos son \emph{Complementos de aritm\'etica y \'algebra}, \emph{\'Algebra elemental moderna} (con J. D. Mancill) y \emph{Fundamentos de la teor\'ia de funciones de variable compleja}. Menos conocido es su primer libro \emph{Introducci\'on al an\'alisis matem\'atico}, publicado en Matanzas en 1940. A poco m\'as de 75 a\~nos de su primera –y \'unica– edici\'on, la situaci\'on de la ense\~nanza de la matem\'atica en Cuba, y en el mundo, pasa por momentos dif\'iciles. Este peque\~no texto puede servir de ejemplo a profesores de diferentes niveles de c\'omo el rigor y la did\'actica no tienen por qu\'e estar en desacuerdo.

En la primera parte del art\'iculo se realiza un peque\~no estudio biogr\'afico de Mario O. Gonz\'alez, donde se comentan algunos datos novedosos obtenidos en el Archivo Hist\'orico Provincial de Matanzas y el Archivo Central de la Universidad de La Habana. En el resto se discuten varios aspectos relacionados con \emph{Introducci\'on al an\'alisis matem\'atico}.


\section{Nota biogr\'afica (1913-40)}

Mario Octavio Jos\'e de Jes\'us de la Caridad Gonz\'alez Rodr\'iguez, naci\'o en la ciudad de Matanzas en la madrugada (5:30 AM) del 14 de septiembre de 1913. Hijo de Mario Gonz\'alez Darna, natural de C\'ardenas, y de Margarita Edelmira Rodr\'iguez Llorente, natural de Matanzas. Sus abuelos maternos y su abuela paterna eran tambi\'en matanceros, mientras que su abuelo paterno era originario de Gij\'on, Asturias. Su padre fue profesor de la c\'atedra B (Ingl\'es) del Instituto de Segunda Ense\~nanza de Matanzas, y adem\'as fungi\'o como Director Interino del mismo en la d\'ecada del 30’. Sus primeros a\~nos los vivi\'o en la vivienda situada en la calle Contreras n\'umero 146, y posteriormente su familia se traslad\'o a una casa con direcci\'on Milan\'es 110\footnote{Hoy en d\'ia la casa de Milan\'es 110 todav\'ia conserva la fachada original.}; ambas situadas a menos de 100m del Instituto de Segunda Ense\~nanza de Matanzas.
 
En los fondos del Archivo Hist\'orico Provincial, se conserva el expediente acad\'emico de Mario O. Gonz\'alez, como alumno del Instituto de Segunda Ense\~nanza de Matanzas. En este, se declara que se matricul\'o para iniciar sus estudios de bachillerato en el a\~no lectivo 1926-27; y los culmin\'o con la obtenci\'on del t\'itulo de bachiller en 1930. El profesor Manuel Labra (1900-82), uno de los m\'as destacados pedagogos cubanos de la primera mitad del siglo \textsc{xx}, fue su profesor de Aritm\'etica, en la que consigui\'o, mediante examen de oposici\'on uno de los premios. Mario O. Gonz\'alez tuvo un desempe\~no sobresaliente en las asignaturas de matem\'atica y f\'isica (\'Algebra, Aritm\'etica, Geometr\'ia y Trigonometr\'ia, F\'isica \textsc{i} y F\'isica \textsc{ii}), en las que obtuvo la m\'axima calificaci\'on, y adem\'as, sus respectivos premios. 
 
La culminaci\'on de la ense\~nanza media y el comienzo de su etapa universitaria, coincidieron con una de las etapas m\'as convulsas de la historia nacional. Entre 1930 y 1933 se produjo un amplio movimiento popular en contra de la dictadura de Gerardo Machado, en la que los estudiantes –especialmente los universitarios–, fueron parte activa. Incluso despu\'es de la fuga del dictador, en agosto de 1933, la agitaci\'on pol\'itica continu\'o reinando en el pa\'is hasta comienzos del a\~no 1937, cuando se «estabiliz\'o» durante el gobierno de Federico Laredo Bru. 

La situaci\'on pol\'itico-social de la naci\'on influy\'o notablemente en la formaci\'on universitaria de Mario O. Gonz\'alez. Su solicitud de ingreso a la Universidad de La Habana se efectu\'o el 16 de septiembre de 1930, con la finalidad de obtener la titulaci\'on de Doctor en Ciencias F\'isico-Matem\'aticas e Ingeniero Electricista. Para el a\~no acad\'emico 1930-31 matricul\'o las asignaturas propias del primer a\~no de la especialidad en Ciencias F\'isico Matem\'aticas, as\'i como Qu\'imica (Inorg\'anica y Anal\'itica) y F\'isica del segundo a\~no. En las materias seleccionadas estaba el An\'alisis Matem\'atico \textsc{i}, que consist\'ia en \'Algebra Superior, pero no tocaba ning\'un tema de \'algebra abstracta.

El inicio del curso estaba previsto para el 3 de diciembre de 1930, pero el movimiento estudiantil en contra de la dictadura de Gerardo Machado mantuvo las aulas universitarias vac\'ias \cite[p. 431]{deArmas}. El dictador trat\'o de disminuir el \'impetu del estudiantado con el cierre de la Universidad de La Habana, por decreto presidencial el 15 de diciembre de 1930, que estuvo vigente hasta mayo de 1932 cuando fue declarado inconstitucional. Sin embargo, el recinto universitario se mantuvo cerrado hasta despu\'es de la ca\'ida del machadato, si bien esta decisi\'on fue tomada por el claustro como medida ante la convulsa situaci\'on pol\'itica y social que se viv\'ia en el pa\'is \cite[p. 453]{deArmas}.

Por las razones anteriores, el primer curso de Mario O. Gonz\'alez fue el correspondiente al a\~no acad\'emico de 1933-34, iniciado el 14 de enero de 1934.  En este, a las asignaturas matriculadas en 1930 se a\~nadieron las del segundo a\~no de Ciencias F\'isico Matem\'aticas. En total fueron 14 materias las cursadas –salvo en Dibujo Natural y F\'isica Superior–, todas sus notas fueron de Sobresaliente; adem\'as, en An\'alisis Matem\'atico \textsc{i} y \textsc{ii}, Trigonometr\'ia y Geometr\'ia obtuvo el Premio de Oro. En su expediente consta que en el inicio de este curso se solicit\'o la matr\'icula para obtener los t\'itulos de Doctor en Ciencias F\'isico Matem\'aticas, F\'isico Qu\'imicas, Ciencias Naturales e Ingeniero Electricista.

En el siguiente a\~no acad\'emico, las asignaturas matriculadas corresponden en su mayor\'ia a las propias del tercer a\~no de estudio de las Ciencias F\'isico-Matem\'aticas. Adem\'as, se encuentran materias m\'as avanzadas como Complementos de An\'alisis Matem\'atico (con Pablo Miquel) y Complementos de Mec\'anica Racional, que fueron de nueva incorporaci\'on al programa de estudios; en ellas se ense\~naron por vez primera de manera sistem\'atica temas como series de Fourier, c\'alculo vectorial y ecuaciones diferenciales \cite{CarlosConcha}. Sin dudas, estos cambios favorecieron su formaci\'on cient\'ifica y crearon un ambiente propicio para un salto cualitativo del desarrollo matem\'atico en Cuba. Tambi\'en curs\'o otras como Agrimensura, Topograf\'ia y Geolog\'ia, m\'as alejadas del perfil matem\'atico.

El fin del machadato no signific\'o el regreso de la estabilidad en la naci\'on, pues cuando parec\'ia que las medidas del Gobierno de los Cien D\'ias lo encaminaban a cumplir, al menos, una parte de los reclamos populares, fue derrocado por el coronel Batista. La dictadura «disimulada» de Batista acarre\'o una nueva ola de protestas, en especial por parte de los estudiantes de la Universidad de La Habana. Por esa raz\'on fue ocupada militarmente a partir de marzo de 1935, y con ello comenz\'o un per\'iodo en el cual la vida universitaria fue sumamente irregular. Aunque el gobierno intent\'o iniciar el curso 1936-37 en mayo de 1936, pero esto no sucedi\'o de manera efectiva hasta marzo de 1937 \cite[pp. 491-493]{deArmas}.

Por el motivo anterior, las notas que aparecen en el expediente de Mario O. Gonz\'alez fueron reportadas en octubre y noviembre de 1937.  Entonces, para fines de ese a\~no, tuvo aprobadas todas las asignaturas correspondientes al plan de estudios del doctorado en Ciencias F\'isico-Matem\'aticas, y solo le faltaban los ejercicios de fin de estudios los que no se efectuaron hasta junio de 1938.

Para obtener el t\'itulo de doctor, adem\'as de aprobar las asignaturas correspondientes, era necesario la discusi\'on de un tema de investigaci\'on; un examen pr\'actico con dos problemas de Matem\'atica y uno de F\'isica; y una breve presentaci\'on con un tema del curr\'iculo escogido por el estudiante. El tribunal formado para el ejercicio de Mario O. Gonz\'alez tuvo como Presidente a Pablo Miquel, como Vocero a Enrique Badell y como Secretario a Rafael Fiterre, el que se reuni\'o el 27 de junio de 1938, para escuchar su tesis «Algunos nuevos tipos de ecuaciones diferenciales invariantes en ciertas transformaciones infinitesimales»; y nuevamente el d\'ia siguiente para el examen pr\'actico, y evaluar una breve conferencia sobre el concepto de derivada. En todos los casos, el tribunal evalu\'o con la m\'axima calificaci\'on los ejercicios, por lo que el t\'itulo de Doctor en Ciencias-F\'isico Matem\'aticas a nombre de Mario O. Gonz\'alez se expidi\'o en octubre de 1938.

Despu\'es de graduarse, imparti\'o clases de Aritm\'etica en el Instituto de Segunda Ense\~nanza de Matanzas y en el de La Habana. En 1939 gan\'o una beca Guggenheim, y el 23 de agosto\footnote{Como consta en una carta firmada por Mario O. Gonz\'alez dirigida al director del Instituto de Segunda Ense\~nanza de Matanzas. Dicha carta se encuentra en expediente laboral de Mario O. Gonz\'alez que se conserva en el Archivo Hist\'orico Provincial de Matanzas.} de ese a\~no parti\'o hacia los Estados Unidos. All\'i realiz\'o estudios de posgrado en el MIT \cite{CarlosConcha}, y al prorrogarse la beca por un a\~no, continu\'o su formaci\'on en la Universidad de Princeton; por lo que su regreso definitivo a Cuba no ocurri\'o hasta bien entrado el a\~no 1941. 

Es en agosto o septiembre de 1940 que se publica, por la Editorial Barani, su libro \emph{Introducci\'on al an\'alisis matem\'atico}. En este momento es becario Guggenheim, por lo que podemos conjeturar que la mayor parte del libro estuvo hecha antes de su viaje a Estados Unidos. Uno de los elementos a favor de esta suposici\'on es que uno de los textos m\'as citados en \emph{Introducci\'on\ldots} –\emph{Theory of functions of a real variable} de E. W. Hobson–, se encuentra en la Biblioteca Central de la Universidad de La Habana desde 1937. Otro de sus referentes fue \emph{Elementos de an\'alisis algebraico} de Rey Pastor, y en la Biblioteca Central se encontraban las ediciones de 1922 y 1939. Aunque no hay evidencias del momento de entrada a los fondos de la biblioteca, es muy probable que, al menos, la edici\'on de 1922 se encontrara a su alcance.

Sin embargo, hay una serie de elementos que apuntan que la estancia en los Estados Unidos le dio la oportunidad de redondear las notas al final de cada cap\'itulo. Esta conjetura se sostiene porque varios textos mencionados en las notas no se encuentran ahora all\'i, y es muy probable que nunca hayan pertenecido a los fondos de la Biblioteca Central. Es el caso de la traducci\'on de Oscar Zariski de \emph{Was sind und was sollen die Zahlen?} o \emph{Elementi di aritmetica regionata e di algebra} de Alfredo Capelli.


\section{Razones para un libro}

El propio Mario O. Gonz\'alez establece que el motivo fundamental que lo movi\'o a escribir \emph{Introducci\'on al an\'alisis matem\'atico} fue la existencia «en nuestra ense\~nanza matem\'atica de algunas discontinuidades entre los estudios secundarios y los universitarios» \cite{MOG40}. 

Seg\'un \'el, para lograr una completa comprensi\'on del an\'alisis matem\'atico, se necesita dominar previamente ciertos temas que no se abordan de forma correcta, o no se tratan en los cursos de matem\'atica elemental. Entre estos destaca la teor\'ia del n\'umero irracional, que es «ineludible antecedente de la teor\'ia de los l\'imites».

En el pr\'ologo de \emph{Complementos de aritm\'etica y \'algebra}, Mario O. Gonz\'alez resalt\'o nuevamente la necesidad del estudio de la teor\'ia del n\'umero irracional: «Puede considerarse la teor\'ia de los l\'imites como la l\'inea fronteriza que separa la matem\'atica elemental de la superior. Pero esta frontera es para algunos, muralla insalvable que les veda el dominio del C\'alculo infinitesimal y otras teor\'ias del An\'alisis. Esta dificultad proviene a nuestro entender, en el aspecto conceptual, del desconocimiento o confusi\'on en que se ha mantenido a los alumnos sobre el important\'isimo concepto de n\'umero irracional» \cite{MOG52}.

Nuevamente, se expone la necesidad de sortear con habilidad pedag\'ogica las dificultades que presenta la ense\~nanza de los n\'umeros irracionales y no evadir el tema. La exposici\'on de estas ideas en \emph{Complementos\ldots} es m\'as reducida que la hecha en \emph{Introducci\'on\ldots}. Desde el punto de vista pr\'actico hay tambi\'en deficiencias, las que ser\'an nuevamente expuestas en \cite{MOG40}: «en nuestra ense\~nanza se elude sistem\'aticamente el c\'alculo aproximado y los m\'etodos de computaci\'on gr\'afica y mec\'anica, y se cultiva el h\'abito de la precisi\'on indeterminada en vez de estudiar los m\'etodos para averiguar en cada caso la aproximaci\'on necesaria y los procedimientos que permitan alcanzarla.» 

En las palabras del autor es sencillo notar que, sobre todas las cosas, prima un incentivo did\'actico m\'as que te\'orico. En el mismo sentido se destaca una separaci\'on notable entre los estudios elementales y los superiores, lo que da muestra de los problemas de la ense\~nanza cubana de la \'epoca. Las opiniones de Mario O. Gonz\'alez sobre las deficiencias de la matem\'atica que se ense\~naba en el bachillerato no provienen solo de su experiencia como alumno, sino tambi\'en como profesor. En 1933 se le abri\'o un expediente laboral en el Instituto de Matanzas, como profesor titular de la c\'atedra de Agricultura. Sin embargo, en el a\~no acad\'emico 1937-38 recibi\'o la autorizaci\'on para impartir Aritm\'etica, quiz\'as gracias a que en ese entonces el director del Instituto era el profesor Manuel Labra. En ese propio curso sirvi\'o como docente de la misma materia en el Instituto No.1 de La Habana.

En los programas de matem\'atica vigentes alrededor de 1940 para los bachilleratos, hab\'ia amplia preponderancia de los temas de geometr\'ia plana y trigonometr\'ia. Los aspectos aritm\'eticos se reduc\'ian a las propiedades b\'asicas de proporcionalidad, divisibilidad de n\'umeros enteros, y algunos temas de matem\'atica financiera. As\'i lo demuestra el contenido de los libros de texto utilizados en el primer per\'iodo republicano \cite{Fiterre49}, \cite{Fiterre52}, \cite{Fiterre55}, \cite{Paz}.

Sin embargo, no debe dejar de buscarse motivaciones dentro de la propia Universidad. A partir del curso de 1912-13 en que Pablo Miquel se convirti\'o en profesor de An\'alisis Matem\'atico en la Universidad de La Habana, el nivel de los cursos creci\'o tanto en rigor como en claridad. Lo anterior se puede evidenciar a trav\'es del texto b\'asico empleado en dichos cursos: \emph{Elementos de \'algebra Superior}, del profesor Miquel. 

La obra de Miquel constituye un hito de la cultura matem\'atica cubana, pues en \'el se expon\'ian los elementos del c\'alculo diferencial e integral con mucho m\'as rigor que en cualquier otro documento elaborado o utilizado en la Isla. Adem\'as, para su momento fue un libro al mismo nivel de actualizaci\'on de sus similares europeos \cite{ConchaCarlos}. Sin embargo, el texto de Miquel ten\'ia una deficiencia: la definici\'on de n\'umero real. Esta cuesti\'on le dificult\'o el tratamiento riguroso, por ejemplo, del teorema de Bolzano sobre los valores intermedios \cite{ConchaCarlos}. 

En \emph{Introducci\'on\ldots} no se menciona de manera expl\'icita el texto de Pablo Miquel, ni la forma en que se trataba el tema de los n\'umeros reales en la Universidad de La Habana. Sin embargo, creemos que el hecho de no existir un tratamiento riguroso de este tema en el texto b\'asico para el curso de c\'alculo diferencial e integral, constituy\'o una motivaci\'on fundamental para que Mario O. Gonz\'alez publicara su peque\~no libro. En buena lid se pod\'ia utilizar alg\'un t\'itulo extranjero, especialmente cuando \emph{Grundlagen der Analysis}\footnote{Fundamentos de an\'alisis}(1929) de Edmund Landau, sigue un gui\'on similar al de \emph{Introducci\'on\ldots}, y constituye un modelo de rigor a seguir para, a partir de los axiomas de los n\'umeros naturales, construir cada uno de los diferentes dominios num\'ericos. Sin embargo, a causa de su rigor extremo no fue bien acogido, ni por los estudiantes, ni por el resto de la comunidad matem\'atica. Por a\~nos este libro ha servido de ejemplo de c\'omo el rigor absoluto puede da\~nar la matem\'atica.

El objetivo de Landau fue desarrollar de forma \textbf{completamente rigurosa} cada sistema num\'erico a partir de los axiomas de los n\'umeros naturales. Siguiendo una idea similar al m\'etodo gen\'etico, construye $\mathbb{Z}$, $\mathbb{Q}$, $\mathbb{R}$ y $\mathbb{C}$. Si bien en el original en alem\'an Landau «tutea» al lector para intentar una comunicaci\'on m\'as cercana, no lo logra. El modo en que est\'a escrito es, «dada la sencillez del contenido, de un estilo despiadadamente similar a un telegrama», con el que estima que, luego de las cinco primeras p\'aginas «de contenido realmente abstracto» la lectura es sencilla.  Tal vez, la mejor opini\'on sobre los\emph{ Grundlagen\ldots} se debe a J. Stillwell: «En mi opini\'on, el problema del libro de Landau no es el rigor (aunque es excesivo) sino la ausencia de notas hist\'oricas, ejemplos y comentarios. Adem\'as, est\'a el hecho de que luego de construir los n\'umeros reales, \'el hace \textbf{absolutamente nada}\footnote{Resaltado del autor.} con ellos.» \cite{Stillwell}.

En contraposici\'on, en el pr\'ologo de \cite{MOG40} aparece: «Muy lejos de nuestro \'animo est\'a pretender alcanzar el rigor absoluto, que no consentir\'ia un libro did\'actico como este.» Este comentario muestra que Mario O. Gonz\'alez est\'a dispuesto a sacrificar un poco de rigor con tal de hacer su obra asequible a un p\'ublico m\'as amplio. Es decir, ¿de qu\'e vale un texto sumamente riguroso que no sea le\'ido? Por ello se justifica plenamente la necesidad de preparar un texto para explicar la teor\'ia del n\'umero irracional a los estudiantes cubanos, interesados en ampliar sus conocimientos matem\'aticos.

Por otro lado, en \cite{Stillwell} se afirma que luego del texto de Landau muy poco o casi nada se escribi\'o sobre la teor\'ia del n\'umero irracional. Esto muestra que \emph{Introducci\'on\ldots} constituye una obra que, a pesar de tratar un tema relativamente b\'asico del An\'alisis, explora una zona poco tratada en la \'epoca: un acercamiento did\'actico a los n\'umeros reales.
Adem\'as de exponer la teor\'ia de n\'umero real, Mario O. Gonz\'alez espera que su breve libro sirva a los profesores de nivel secundario, dado la inclusi\'on de numerosos comentarios hist\'oricos. De esta forma, la utilizaci\'on de la historia como recurso did\'actico por los profesores de matem\'atica en Cuba puede encontrar un precedente en la obra de esta eminente figura.

\section{Interioridades del texto}

El primero de los cap\'itulos de \emph{Introducci\'on al An\'alisis Matem\'atico} est\'a dedicado a los n\'umeros naturales. Este comienza con un ep\'igrafe de F\'elix Klein: «Las propiedades b\'asicas de los n\'umeros y las operaciones que con ellos se efect\'uen pueden reducirse a las propiedades generales de los conjuntos y a las relaciones abstractas que entre ellos existen». Haciendo honor a estas palabras, 17 de las 33 secciones est\'an dedicadas a estudiar propiedades b\'asicas de la teor\'ia de conjuntos que servir\'an para deducir las reglas operativas de los n\'umeros naturales. En este sentido, la obra de Mario O. Gonz\'alez est\'a claramente influida por las ideas de Dedekind en \emph{Was sind und was sollen die Zahlen?}.

Como libro did\'actico no se consideran los axiomas de la teor\'ia de conjuntos, por lo que se deja a un nivel intuitivo. Sin embargo, no se da pie a que el tratamiento \emph{na\"{\i}ve} d\'e lugar a las antinomias, pues se deja claro que no se puede «subordinar la definici\'on de conjunto a la de sus elementos y viceversa».  Esto muestra que Mario O. Gonz\'alez no solo se encuentra al tanto de los problemas de los fundamentos de la teor\'ia de conjuntos, sino que los comprende y expone con claridad.

Son varios los puntos de contacto que tiene \emph{Introducci\'on\ldots} con \emph{Was sind\ldots}, pues all\'i se utiliza el concepto de conjunto simplemente infinito y se demuestra el principio de inducci\'on matem\'atica. Pero tambi\'en existen diferencias notables, pues Mario O. Gonz\'alez define conjunto finito y utiliza el principio de buen ordenamiento como base de su teor\'ia. Esto le permite simplificar todo el aparataje te\'orico que Dedekind introdujo en su monograf\'ia de 1888.
 
Para Mario O. Gonz\'alez el n\'umero natural es primero cardinal, a diferencia de Dedekind para quien es antes ordinal. Sin embargo, prevalece la abstracci\'on para definir el concepto de n\'umero natural: «N\'umero natural es, pues, el concepto abstracto nacido espont\'aneamente en nuestra conciencia para representar los conjuntos que son coordinables\footnote{L\'ease equipotentes.} y diferenciarlos de los no coordinables.» \cite[p. 14]{MOG40}

Las operaciones aritm\'eticas entre n\'umeros naturales se definen a partir de operaciones entre conjuntos. La suma se define a partir de la uni\'on, y la multiplicaci\'on a partir del producto cartesiano. Aqu\'i vale se\~nalar que se utiliza el t\'ermino \emph{composici\'on de conjuntos} en vez de producto cartesiano. Adem\'as, este m\'etodo oscurece un poco la definici\'on de las operaciones aritm\'eticas respecto a Dedekind, y se debe a que Mario O. Gonz\'alez prescinde de la operaci\'on \emph{sucesor} utilizada en \emph{Was sind\ldots} para definir la suma, multiplicaci\'on y elevaci\'on a potencias naturales.

Como colof\'on del primer cap\'itulo, aparece una secci\'on de \textbf{Notas}, en las que se exponen varios comentarios hist\'oricos. Aqu\'i se hace un breve recorrido por el concepto de n\'umero natural, en la que se destaca el papel de Dedekind y su \emph{Was sind\ldots}, la que considera como la m\'as aceptada. Tambi\'en expone los axiomas de Peano para $\mathbb{N}$ y los de Hilbert para $\mathbb{R}$. Estas \textbf{Notas} muestran la amplia cultura matem\'atica, el especial inter\'es por la historia de esta ciencia y los problemas de sus fundamentos. 

En el segundo cap\'itulo se tratan los n\'umeros racionales positivos, los que «de acuerdo a su origen hist\'orico y teniendo presente atendibles razones de orden did\'actico se establece el concepto de n\'umero racional apoy\'andonos en la teor\'ia de la medida. Ello tiene otra ventaja y es que permite fijar con precisi\'on los conceptos de magnitud y cantidad, los cuales entre nosotros suelen tratarse de manera bastante err\'onea y confundirse muchas veces.»\cite[p. 41]{MOG40}.

Se puede comprobar aqu\'i como Mario O. Gonz\'alez no separa la idea abstracta de n\'umero de su origen pr\'actico. Es decir, no se desliga la matem\'atica de sus ra\'ices en la actividad humana, pues el origen hist\'orico de los n\'umeros racionales positivos est\'a en la acci\'on de medir. Por ello se han utilizado para representar las cantidades de las m\'as diversas magnitudes, y entre estas \'ultimas, la m\'as sencilla es \emph{longitud}; la que tambi\'en se ha utilizado como representaci\'on «natural» de los n\'umeros racionales. Estas ideas le permiten introducir de forma natural el c\'alculo con n\'umeros decimales aproximados. Adem\'as, no se sigue el m\'etodo gen\'etico estrictamente, pues los n\'umeros racionales se definen a partir de $\mathbb{N}$, no se construyen. 

Es en este cap\'itulo donde resalta sobremanera el hecho de que Mario O. Gonz\'alez no desliga el concepto de n\'umero de su utilidad pr\'actica, y sobre todo de la profunda simbiosis existente entre el «mundo matem\'atico» y el «mundo real». Por ejemplo, veamos su comentario respecto a la posibilidad de definir de otra manera las operaciones aritm\'eticas entre los n\'umeros racionales:
\begin{quote}
 desde el punto de vista puramente l\'ogico podr\'ian ser adoptadas otras definiciones que, presentando esa simetr\'ia formal, dejaran a salvo el principio de permanencia de Hankel. Mas, sin duda, el olvido que tales definiciones supondr\'ian del significado concreto de los n\'umeros, har\'ia perder a la teor\'ia que se desarrollase sobre ellas todo inter\'es cient\'ifico; con raras excepciones, las doctrinas que se han edificado desde\~nando el contacto con la realidad objetiva, han tardado poco en quedar anquilosadas.

Las definiciones usuales no responden, pues, a una necesidad l\'ogica fatalmente inevitable; s\'i a una obligatoriedad pr\'actica imperativa. \cite[p. 59]{MOG40}.
\end{quote}

No se le debe dar a estas palabras un sentido puramente literal, pues pudiera cometerse el error de malinterpretarlas. Es claro que en este caso no se refiere al desarrollo de nuevas teor\'ias a partir de la generalizaci\'on de ideas y conceptos ya conocidos. Tambi\'en es claro que, si se sacan de contexto estas palabras, la posici\'on de Mario O. Gonz\'alez estar\'ia en franca contradicci\'on con la introducci\'on de los n\'umeros complejos o los cuaterniones.

El tercer cap\'itulo est\'a dedicado a la introducci\'on de los n\'umeros negativos, la que se realiza tanto a partir de consideraciones aritm\'eticas como pr\'acticas. Desde el punto de vista aritm\'etico los n\'umeros negativos son necesarios para dar sentido a la resta entre dos n\'umero racionales arbitrarios. En sentido pr\'actico, los n\'umeros negativos aparecen al considerarse magnitudes que pueden medirse en dos sentidos, como es el caso de la relaci\'on entre deudor y acreedor.

Pero sin duda, el cap\'itulo m\'as interesante es el cuarto, porque es el dedicado al estudio de los n\'umeros reales que, en esencia, es un \emph{leit motiv} de \emph{Introducci\'on\ldots}. El cap\'itulo se inicia con un ep\'igrafe de J. Tannery: «Yo he adoptado la idea fundamental de Dedekind, que me parece ilumina profundamente la naturaleza del n\'umero irracional.». Con ella se anuncia el camino que se seguir\'a para la construcci\'on de los n\'umeros reales: las cortaduras de Dedekind.
 
La primera secci\'on es un peque\~no comentario hist\'orico en el que se demuestra la inconmensurabilidad de la diagonal del cuadrado de lado 1, como «descubrimiento» de los n\'umeros irracionales. Seguidamente se muestra que para resolver el problema de logaritmaci\'on y elevaci\'on a potencia racional, se necesitan los n\'umeros irracionales. De esta forma Mario O. Gonz\'alez muestra que la introducci\'on de los irracionales parte de una necesidad pr\'actica.

De acuerdo al ep\'igrafe que da inicio al cap\'itulo, se escoge el m\'etodo de las cortaduras de Dedekind. Esta elecci\'on se debe a que de las varias v\'ias para construir $\mathbb{R}$, las cortaduras es la m\'as sencilla. Sin embargo, el punto d\'ebil del m\'etodo de Dedekind es lo enrevesado que resulta definir las operaciones aritm\'eticas fundamentales. Luego, por razones did\'acticas, Mario O. Gonz\'alez escoge otro camino para definir las operaciones aritm\'eticas en el conjunto de los n\'umeros reales: «Al exponer las operaciones aritm\'eticas nos apartamos del m\'etodo de Dedekind porque estimamos m\'as claro y tambi\'en m\'as pr\'actico ense\~nar a ejecutarlas utilizando las expresiones decimales de los n\'umeros reales» \cite[p. 81]{MOG40}. De esta manera se comprueba que el objetivo de Mario O. Gonz\'alez no es preparar un texto con la construcci\'on de los n\'umeros reales, sino hacerlo de forma did\'actica y sencilla.

Es necesario resaltar la profunda comprensi\'on de las ideas de Dedekind que logra Mario O. Gonz\'alez. Muestra de ello es un breve comentario que aparece en el ac\'apite dedicado al concepto de n\'umero real:
\begin{quote}
 Algunos autores consideran la cortadura en s\'i como el n\'umero irracional, de suerte que para ellos hablar de n\'umero irracional o de una cortadura sin elemento de separaci\'on en el campo de los n\'umeros racionales, es exactamente la misma cosa. Un n\'umero irracional   no ser\'ia, pues, sino un s\'imbolo con el cual se representar\'ian las dos clases de una cortadura del tipo mencionado.

Aunque esta manera de ver la cuesti\'on sea perfectamente leg\'itima, nos parece, sin embargo, m\'as de acuerdo con el pensamiento original de Dedekind [\ldots] considerar la cortadura $(A,B)$ como representaci\'on de un ente abstracto (el n\'umero irracional), creado para servir de elemento de separaci\'on de las dos clases, y no inversamente. \cite[p. 84]{MOG40}
\end{quote}
 Esto demuestra que Mario O. Gonz\'alez no se contenta con la reproducci\'on de las ideas de Dedekind, sino que se preocupa por comprender las ideas fundacionales desarrolladas por el matem\'atico germano.

El basamento te\'orico escogido por Mario O. Gonz\'alez para simplificar el tratamiento de las operaciones aritm\'eticas en $\mathbb{R}$ es el concepto de \emph{clases contiguas}\footnote{Es necesario se\~nalar que Mario O. Gonz\'alez demuestra rigurosamente que una cortadura define dos clases contiguas, y que dos clases contiguas define una cortadura.}. El origen de este m\'etodo se debe a Alfredo Capelli entre 1897 y 1903. Este m\'etodo evita la utilizaci\'on de todos los elementos de la cortadura, basta considerarse una sucesi\'on en cada clase con la condici\'on que la diferencia entre sus t\'erminos sea infinitesimal. Luego, en el lenguaje de las sucesiones es mucho m\'as sencillo lidiar con las operaciones aritm\'eticas.

Las operaciones aritm\'eticas discutidas por Mario O. Gon\-z\'alez son la suma, resta, multiplicaci\'on, divisi\'on, radicaci\'on, potencia con exponente racional y real, y la logaritmaci\'on. Este espectro de operaciones en $\mathbb{R}$ es superior al estudiado en \cite{Landau}; sin embargo, no es un aporte de su persona, porque ya est\'a en \cite{ReyPastor48}\footnote{La edici\'on de 1922 de este texto se encontraba en la Bilblioteca Central de la Universidad de La Habana, porque es muy probable que haya sido consultado por Mario O. Gonz\'alez.}. 

La secci\'on «Algunas propiedades del conjunto de los n\'umeros reales» est\'a dedicada al estudio de propiedades menos operativas de $\mathbb{R}$. En ella se enuncia y demuestra correctamente el teorema de Dedekind referente a que una cortadura en $\mathbb{R}$ determina un \'unico n\'umero real, y viceversa. Adem\'as, se introduce la noci\'on de conjunto numerable y se prueba la numerabilidad de los n\'umeros racionales y la no numerabilidad de los reales contenidos en el $(0,1)$.
 
En las notas al cap\'itulo, se critica el concepto de n\'umero irracional dado por Aurelio Baldor en \emph{Aritm\'etica te\'orico pr\'actica}, que por mucho tiempo fue utilizado en las escuelas cubanas y latinoamericanas. Sin embargo, lo m\'as destacable es el estudio hist\'orico del concepto de n\'umero real. En este bosquejo se discuten diferentes teor\'ias aparecidas en el siglo \textsc{xix} para explicar la naturaleza del n\'umero real, como las m\'as importantes debidas a Weierstrass, Cantor y Dedekind; y otras menos conocidas como las de C. M\'eray. Tambi\'en se tratan otras a\'un menos divulgadas, como son las aportadas por C. Arzel\'a, P. Bachmann y A. Capelli, y adem\'as se incluye un breve comentario de las ideas promovidas por Kronecker. Es curioso que en su breve «historia del n\'umero real», Mario O. Gonz\'alez no menciona a Heine, -quien dio una construcci\'on de los n\'umeros reales en 1872 de forma similar a Cantor-; ni la versi\'on axiom\'atica de Hilbert, aunque es necesario se\~nalar que en el primer cap\'itulo del libro aparecen los axiomas de la aritm\'etica dados por este \'ultimo. En estas notas al cap\'itulo, se destacan dos elementos esenciales de la obra de Mario O. Gonz\'alez, especialmente por ser \emph{Introducci\'on\ldots} una obra de juventud: su amplia cultura matem\'atica y su especial inter\'es por la historia de la ciencia.

Otro de los aspectos sobresalientes de \emph{Introducci\'on\ldots} es su particular \'enfasis en el c\'alculo con n\'umeros aproximados. ¿Por qu\'e se necesita tratar este tema? Seg\'un Mario O. Gonz\'alez porque: «no obstante su gran importancia te\'orica y pr\'actica, a pesar de lo extraordinariamente \'util que es al matem\'atico puro, al f\'isico, al qu\'imico, al agrimensor, al ingeniero, y en general, a todos los que hacen aplicaci\'on de la ciencia matem\'atica, el c\'alculo con n\'umeros aproximados se ha visto hasta ahora excluido sistem\'aticamente  de nuestros programas oficiales.» \cite[p. 133]{MOG40}. Es decir, la necesidad es eminentemente pr\'actica. La otra lectura que se le debe hacer es sobre las deficiencias del sistema educacional cubano de la \'epoca. En un pa\'is donde uno de los objetivos fundamentales de ense\~nar matem\'atica son sus aplicaciones y adem\'as, en un momento donde no exist\'ian a\'un las grandes m\'aquinas de c\'omputo, sus educandos no ten\'ian dentro de su programa la herramienta que les permit\'ia obtener los resultados num\'ericos necesarios.

Los n\'umeros complejos son introducidos de manera similar a la construcci\'on de Hamilton. Las definiciones de las operaciones aritm\'eticas b\'asicas en este dominio son deducidas a partir del principio de permanencia de las leyes formales de Hankel. Este m\'etodo posibilita la deducci\'on de las operaciones en $\mathbb{C}$, a partir de sus similares reales. Adem\'as, de cierta forma se evita la introducci\'on «axiom\'atica» de las operaciones con n\'umeros complejos, para luego verificar que son compatibles con las definidas en $\mathbb{R}$. 

Lo que resulta extra\~no es el motivo por el cual no se da un orden para el sistema de los n\'umeros complejos: «porque en este campo carecen de importancia» afirma Mario O. Gonz\'alez en \cite[p. 155]{MOG40}. Una explicaci\'on razonable para esto es evitar el an\'alisis de la compatibilidad de las operaciones aritm\'eticas con el orden. En el estudio de los n\'umeros complejos que se realiza en \cite{MOG52} no hay ning\'un comentario sobre la no existencia de orden en $\mathbb{C}$.

En el resto del cap\'itulo se estudian la forma binomial y la representaci\'on geom\'etrica de los n\'umeros complejos, as\'i como la interpretaci\'on geom\'etrica de la suma y la multiplicaci\'on. Se introduce la forma trigonom\'etrica, la f\'ormula de Moivre y la extracci\'on de ra\'ices. El cap\'itulo –y por ende el libro– con el estudio de la elevaci\'on de un n\'umero complejo a potencias no enteras, y la definici\'on del logaritmo complejo.

En las notas al cap\'itulo se hace un breve recorrido por la historia de los n\'umeros complejos, en \'el se resalta c\'omo demor\'o su aceptaci\'on por la comunidad matem\'atica. Adem\'as, se realiza una peque\~na introducci\'on a los cuaterniones.

\section{A modo de conclusi\'on}

La trascendencia de \emph{Introducci\'on al An\'alisis Matem\'atico} debe buscarse en varios sentidos, y no porque sea el primer texto de Mario O. Gonz\'alez. Este libro es, casi seguramente, el primero impreso fuera de La Habana en lidiar con temas de matem\'atica. Adem\'as, aborda un asunto que, aunque elemental, no es nada sencillo, y lo hace con mucho didactismo. La introducci\'on rigurosa de los n\'umeros reales nunca ha sido un aspecto trivial del an\'alisis, dada su profunda relaci\'on con los fundamentos de la matem\'atica. 

Quiz\'as, fue tambi\'en el primer texto cubano de matem\'atica en trascender m\'as all\'a de las fronteras nacionales. As\'i lo demuestra el siguiente comentario de Rey Pastor sobre \emph{Introducci\'on\ldots}: «ensayo de cr\'itica, en algunos puntos atinada, de las diversas teor\'ias del n\'umero» \cite[p. 6]{ReyPastor48}. Esta opini\'on evolucion\'o, y en su monumental obra \emph{An\'alisis Matem\'atico} Tomo I comenta que  es «un ensayo de cr\'itica  de las diversas teor\'ias del n\'umero. Con exposici\'on preferente de la teor\'ia cardinal, conteniendo en cada paso adecuada rese\~na hist\'orica y adem\'as un cap\'itulo sobre n\'umeros aproximados de gran valor did\'actico». 

No se realizaron otras ediciones porque casi todos los contenidos tratados pasaron a formar parte de \cite{MOG52}. Por ejemplo, el contenido del cap\'itulo \textsc{i} de \emph{Introducci\'on\ldots} pasaron al cap\'itulo \textsc{i} de \emph{Complementos\ldots}; mientras que el estudio de los n\'umeros reales, el c\'alculo con n\'umeros aproximados y los n\'umeros complejos, se encuentran en los cap\'itulos \textsc{iii} y \textsc{vi} de \emph{Complementos\ldots}, respectivamente.

Por otro lado, el texto de Mario O. Gonz\'alez se inscribe en lo que puede ser llamada la tradici\'on cubana de textos de an\'alisis matem\'atico, inaugurada con \emph{Elementos de \'algebra superior} de Pablo Miquel, quien la continu\'o con sus \emph{Curso de c\'alculo diferencial} (1941) y \emph{Curso de c\'alculo integral} (1942). Adem\'as de iniciar otro aspecto importante de la matem\'atica cubana: el inter\'es por la historia de la matem\'atica en s\'i, y su utilidad en la ense\~nanza.

------------------------------------------------
\phantomsection
\section*{Fuentes documentales} 

\begin{itemize}
\item Archivo Central de la Universidad de La Habana: Expediente acad\'emico de Mario O. Gonz\'alez Rodr\'iguez. 
\item Archivo Hist\'orico Provincial de Matanzas: Expediente acad\'emico de Mario O. Gonz\'alez Rodr\'iguez, Instituto de Segunda Ense\~nanza de Matanzas.
\end{itemize}

%
%
%
\bibliographystyle{plain}
\bibliography{Mario_O_Gonzalez}

\begin{thebibliography}{10}

\bibitem{BatardVillegas}
L.~F. Batard~Mart\'inez and P.~J. Villegas~Aguilar.
\newblock {\em Las ciencias exactas y naturales en Cuba}.
\newblock Editorial Cient\'ifico-T\'ecnica, La Habana, 2010.

\bibitem{deArmas}
R.~de~Armas et~al.
\newblock {\em Historia de la Universidad de La Habana. 1930-1978}.
\newblock Editorial de Ciencias Sociales, La Habana, 1984.

\bibitem{Fiterre49}
I.~Fiterre~Riveras.
\newblock {\em Matem\'atica. Segundo curso. Geometr\'ia}.
\newblock Editorial Selecta, La Habana, 4 edition, 1949.

\bibitem{Fiterre52}
I.~Fiterre~Riveras.
\newblock {\em Matem\'atica. Tercer curso. Geometr\'ia y nociones de
  trigonometr\'ia}.
\newblock Editorial Selecta, La Habana, 5 edition, 1952.

\bibitem{Fiterre55}
I.~Fiterre~Riveras.
\newblock {\em Matem\'atica. Cuarto curso. Geometr\'ia}.
\newblock Editorial Selecta, La Habana, 5 edition, 1955.

\bibitem{MOG40}
M.~O. Gonz\'alez.
\newblock {\em Introducci\'on al an\'alisis matem\'atico}.
\newblock Editorial Barani, Matanzas, 1940.

\bibitem{MOG52}
M.~O. Gonz\'alez.
\newblock {\em Complementos de an\'alisis y \'algebra}.
\newblock Editorial Selecta, La Habana, 8 edition, 1952.

\bibitem{Landau}
E.~Landau.
\newblock {\em Foundations of analysis}.
\newblock Chelsea Publishing Co., New York, 1966.

\bibitem{Miyares}
A.~Miyares and J.~M. Escalona.
\newblock {\em Matem\'atica segundo curso. Geometr\'ia}.
\newblock Empresa Consolidada de Artes Gr\'aficas, La Habana, 3 edition, 1962.

\bibitem{Paz}
A.~Paz~Sord\'ia.
\newblock {\em Matem\'atica cuarto curso. Geometr\'ia}.
\newblock Casa Lori\'e, La Habana, 1951.

\bibitem{ReyPastor48}
J.~Rey~Pastor.
\newblock {\em Elementos de an\'alisis algebraico}.
\newblock Soc. de Resp. Ltda., Buenos Aires, 3 edition, 1948.

\bibitem{ReyPastor69}
J.~Rey~Pastor et~al.
\newblock {\em An\'lisis matem\'atico}, volume~1.
\newblock Editorial Kapelusz, Buenos Aires, 8 edition, 1969.

\bibitem{CarlosConcha}
C.~S\'anchez~Fern\'andez and C.~Vald\'es~Castro.
\newblock Emergencia de una cultura matem\'atica en {C}uba.
\newblock {\em Ciencias Matem\'aticas}, 27(2):3--16, 2013.

\bibitem{Stillwell}
J.~Stillwell.
\newblock {\em The real numbers}.
\newblock Springer Verlag, Heidelberg, 2013.

\bibitem{ConchaCarlos}
C.~Vald\'es~Castro and C.~S\'anchez~Fern\'andez.
\newblock El texto \emph{Elementos de \'Algebra Superior} de {P}ablo {M}iquel:
  un salto cualitativo en la instrucci\' n matem\'atica en {C}uba.
\newblock {\em Ciencias Matem\'aticas}, 20(1):54--66, 2002.

\end{thebibliography}
\nocite{BatardVillegas} \nocite{Landau} \nocite{Miyares} \nocite{ReyPastor69}

\end{document}